\documentclass[a4paper,12pt,onecolumn,notitlepage]{article}
\usepackage[pdftex,final]{graphicx}
\usepackage[pdfpagemode=None,bookmarksopen=false,colorlinks=true,linkcolor=blue,citecolor=blue]{hyperref}
\usepackage[all]{xy}
\usepackage{amsmath}
\usepackage{amssymb}
\usepackage{sectsty}
\usepackage{bibspacing}
\usepackage{rpkid}
\long\def\symbolfootnote[#1]#2{\begingroup%
\def\thefootnote{\fnsymbol{footnote}}\footnote[#1]{#2}\endgroup}
\allsectionsfont{\normalsize}
\begin{document}
\thispagestyle{empty}
\numberwithin{equation}{section}
\begin{center}
\noindent
{\Large \textbf{On principles of inductive inference\symbolfootnote[1]{\scriptsize To be published in: Goyal P. (ed.),  \textit{Proceedings of the 31th International Workshop on Bayesian Inference and Maximum Entropy Methods in Science and Engineering, Waterloo, July 10-15, 2011}, AIP Conf. Proc., Springer, Berlin.}}}\\
\ \\
{Ryszard Pawe{\l} Kostecki}\\
{\small\ \\}{\small
\textit{Institute of Theoretical Physics, University of Warsaw, Ho\.{z}a 69, 00-681 Warszawa, Poland}}{\small\\
\ \\
January 31, 2012}
\end{center}
\begin{abstract}{\small\noindent
We propose an \textit{intersubjective epistemic} approach to foundations of probability theory and statistical inference, based on relative entropy and category theory, and aimed to bypass the mathematical and conceptual problems of existing foundational approaches.}
\end{abstract}

\section{Standard frameworks} 

The notion \textit{inference} means `logical reasoning'. The \textit{deductive inference} specifies premises by the valuations of sentences in \textit{truth} values, and provides an inference procedure which is considered to lead to \textit{certain} conclusion on the base of given premises. The \textit{inductive inference} specifies premises by the valuations of sentences in \textit{possible} (\textit{plausible}) values and provides an inference procedure which is considered to lead to \textit{most possible} (\textit{most plausible}) conclusions on the base of given premises. From the mathematical perspective, the difference between deductive and inductive inference lays not in the form of logical valuations (these can be the same in both methods), but in the procedure of specifying conclusions on the base of premises. The conclusions of multiple application of deductive inference to the sequence of sets of premises depend \textit{in principle} on all elements of all these sets, while the conclusions of the multiple application of inductive inference to the sequence of sets of premises depend \textit{in principle} only on some elements of some of these sets. For  this reason, the premises of inductive inference are also called \textit{evidence}. An example of inductive inference procedure is any statistical reasoning based on probabilities. The evidence can be provided, for example, by particular quantities with units interpreted as `experimental data' \textit{together with} a particular choice of a method which incorporates these `data' into statistical inference. Any choice of such method \textit{defines} the evidential meaning of the `data', and is a crucial element of the inference procedure. A standard example of such method is to ignore everything what is known about a sequence of numbers associated with a single abstract quality (such as a ``position''), leaving only the value of arithmetic average and the value of a fluctuation around this average as a subject of comparison (e.g., by identification) with the mean and variance parameters of the gaussian probabilistic model.

According to frequentist interpretation (by Ellis \cite{Ellis:1843}, Venn \cite{Venn:1866}, Fisher \cite{Fisher:1922}, von~Mises \cite{vonMises:1919,vonMises:1928}, Neyman \cite{Neyman:1937}, and others) probabilities can be given meaning only as relative frequencies of some experimental outcomes in some asymptotic limit. This interpretation was very influential in the last 160 years, \textit{but} so far \textit{none} mathematically strict and logically sound formulation of this interpretation exists (see, e.g., \cite{Jeffrey:1977,vanLambalgen:1987,Hajek:1997}). The separation of a formalism of inductive inference into `probability theory' and `statistics' is also a consequence of frequentist interpretation, which forbids consideration of probability (understood as relative frequency) as a subject of change based on change of evidence. \textit{Thus, without frequentism, there is no reason for keeping the division of probabilistic and statistical part of the framework of inductive inference into two separated theories.} Moreover, the methods of statistical inference used within the frames of the frequentist approach are mainly based on \textit{ad hoc} principles, which are justified by convention, and do not possess mathematically strict and logically sound justification (see e.g. \cite{Pratt:1961,Jaynes:1976,Berger:Berry:1988,Shiryayev:Koshevnik:1989,Jaynes:2003}). This is a consequence of the lack of strict and sound foundations of frequentist interpretation of probability.

Beyond logically and mathematically unjustified frequentist approach (and even less successful \cite{Humphreys:1985,Eagle:2004} propensity interpretation \cite{Peirce:1910,Popper:1959,Giere:1973}), probability theory and statistical inference theory can be considered as two parts (resp., kinematic and dynamic) of a single theory of quantitative inductive inference. The evidences used in this inference need not be restricted to frequencies. Nevertheless, the choice of particular methods of specification of evidences and drawing inferences requires some justification. 

The `subjective' bayesian approach (by Ramsey \cite{Ramsey:1931}, de~Finetti \cite{deFinetti:1931,deFinetti:1937,deFinetti:1970}, Savage \cite{Savage:1954}, and others) allows any kinematics and requires Bayes' rule as dynamics, grounding both in requirement of \textit{personal} consistency of betting. This is conceptually consistent, \textit{but} by definition \textit{lacks} any rules relating the probability assignments (theoretical model construction) with intersubjective knowledge (experimental setup construction, `experimental data'). Thus, it is often accused of arbitrariness. Such accusations \textit{are} justified if they amount to saying that the methods of scientific inquiry \textit{seem to be something more} than individual personal consistency of bets, but \textit{are not} justified if they appeal to (operationally undefined!) notions of `objectivity', `nature', `reality', etc., because \textit{any} theoretical statement is after all an arbitrary mental construct. 

The syntactic approach (by Johnson \cite{Johnson:1921}, Keynes \cite{Keynes:1921}, Carnap \cite{Carnap:1950,Carnap:1952}, and others) amounts to construct probability theory as a predicative calculus in a formal language, \textit{but} it does \textit{not} provide neither any sound justification for the choice of language and calculus nor any definite methods of model construction, different from `subjective' approach (see e.g. \cite{Howson:Urbach:1989,Hajek:2009,Zabell:2009}). This makes syntactic approach foundationally irrelevant.

The `objective' bayesian approach (by early Jeffreys \cite{Jeffreys:1931}, Cox \cite{Cox:1946,Cox:1961}, Jaynes \cite{Jaynes:2003}, Berger \cite{Berger:2006}, and others) allows various mathematical rules of assignment of probabilities (see e.g. \cite{Kass:Wasserman:1996,Jaynes:2003}) and of inference (see e.g. \cite{Cox:1946,Cox:1961,Shore:Johnson:1980,Shore:Johnson:1981,Caticha:2004,Caticha:2007}). It attempts to select the preferred rules by an appeal to some notions of `rational consistency' or `experimental reproducibility',  \textit{but} it \textit{fails} to provide sound conceptual justification for these rules which would be neither subjectively idealistic (personalist) nor ontologically idealistic (frequentist) \cite{Uffink:1995,Uffink:1996,Fienberg:2006,Kadane:2006}. (One of important problems is the lack of justification of the methods that associate particular constraints $F$ in entropic updating rule with some functions of `experimental evidence', see e.g. \cite{Uffink:1995,Uffink:1996}.) Yet, the idea to provide some rules of probabilistic model construction taking into account the role played by experimental evidence in \textit{intersubjective consistency} of inductive inferences seems to be crucial.

The variety of above \textit{conceptual frameworks} for probability theory and statistical inference corresponds to the variety of \textit{mathematical frameworks}. There are four main approaches: by Bayes--Laplace \cite{Bayes:1763,Laplace:1812}, Borel--Kolmogorov \cite{Borel:1909:denom,Kolmogorov:1933}, Whittle \cite{Whittle:1970,Whittle:2000}, and Le~Cam \cite{LeCam:1964,LeCam:1986}. But even more approaches is possible, because probability theory can be built upon two components: evaluational (kinematic) and relational (dynamic), and, apart from selection of one or two of these components, one can provide different mathematical implementations thereof. For example, the evaluational component can be given either by an abstract measure theory on $\sigma$-algebras of subsets of some set, or by an integral theory on vector lattices, while the relational component might be given either by Bayes' rule, or by conditional expectations, or by entropic updating, etc. 

The Borel--Kolmogorov framework \cite{Borel:1909:denom,Kolmogorov:1933} is based on the notions of measure spaces $(\X,\mho(\X),\mu)$ and probabilistic models 
\[
	\M(\X,\mho(\X),\mu)\subseteq L_1(\X,\mho(\X),\mu)^+.
\]
Building upon measure-theoretic integration theory, this framework is, from scratch, equipped with kinematic (evaluational) prescriptions, \textit{but} has no \textit{generic} notion of conditional updating of probabilities. (The reason of it is an associated, but by no means necessary, frequentist interpretation, which claims identification of probabilities with frequencies. This forbids `updating' probabilities because it would mean updating the frequencies.) There are three facts to observe here. First, many different measure spaces $(\X,\mho(\X),\mu)$ lead to $L_1(\X,\mho(\X),\mu)$ spaces that are all isometrically isomorphic to the same abstract $L_1(\mho)$ space, where $\mho$ is camDcb-algebra (countably additive, Dedekind complete, boolean, allowing at least one strictly positive semi-finite measure), constructed by \cite{Fremlin:2000}
\[
	\mho:=\mho(\X)/\{A\in\mho(\X)\mid\mu(A)=0\}.
\]
\textit{Thus, only $L_1(\mho)$ is necessary for defining probabilistic models}. But, given any camDcb-algebra $\mho$, the association of $L_1(\mho)$ (and any other $L_p(\mho)$) to $\mho$ is functorial \cite{Fremlin:2000}, and no appeal to representations in terms of measure space is ever required. Second, by the Loomis--Sikorski theorem \cite{Loomis:1947,Sikorski:1948}, each camDcb-algebra $\mho$ can be represented as a measure space $(\X,\mho(\X),\mu)$, given the choice of measure $\mu$ on $\mho$. However, there are many different measure spaces that lead to the same algebra $\mho$ \cite{Segal:1951}. \textit{Thus, using the measure spaces $(\X,\mho(\X),\mu)$ instead of camDcb-algebras $\mho$ is ambiguous.} Finally, as observed by Le~Cam \cite{LeCam:1964,LeCam:1986} and Whittle \cite{Whittle:1970,Whittle:2000}, probabilistic description in terms of measures $\mu$ on $(\X,\mho(\X))$ can be completely replaced by the description in terms of integrals $\omega$ on vector lattices $\stone(\mho)$. (For every camDcb-algebra $\mho$ there exists a canonically associated vector lattice $\stone(\mho)$ of characteristic functions on the set of boolean homomorphisms $\mho\ra\ZZ_2$.) \textit{Thus, one can deal exclusively with positive finite integrals over commutative algebras $\mho$ instead of measures on $\mho(\X)$}. The normalised integrals are just expectation functionals, and probability of $a\in\mho$ is recovered by evaluation $p(a):=\omega(\chr_a)$ on characteristic function $\chr_a\in\stone(\mho)$. 

On the other hand, the Bayes--Laplace framework \cite{Bayes:1763,Laplace:1812} is based on finitely additive boolean algebras $\B$ and conditional probabilities $p(x|\theta):\B\times\B\ra[0,1]$. It is equipped from scratch with dynamical (relational) prescription, given by Bayes' rule 
\[
	p(x|\theta)\mapsto p_{\mathrm{new}}(x|\theta):=p(x|\theta)\frac{p(b|x\land \theta)}{p(b|\theta)},
\]
\textit{but} it provides no \textit{generic} notion of probabilistic expectation over a continuous (countably additive) domain of infinite sets. Bayes' rule defines a concrete method of providing statistical inferences. Thus, \textit{statistical inference can be understood as a dynamical component of probability theory}. As noticed by Jaynes \cite{Jaynes:2003}, Bayes' rule and all Bayes--Laplace framework has precisely the same properties and the same range of validity if the conditional probabilities are evaluated into $[1,\infty]$ instead of $[0,1]$. \textit{Hence, the normalisation of probabilities is not a necessary feature of probabilistic/inferential framework} (it is necessary only in the frequentist interpretation). Moreover, \textit{Bayes' rule is a special case of constrained relative entropic updating} \cite{Caticha:Giffin:2006}  
\[
	p(x|\theta)\mapsto p_{\mathrm{new}}(x|\theta):=\arginf_{q\in\M}\left\{\int p\log(p/q)+F(q)\right\},
\] 
for $p,q$ in a probabilistic model $\M$, $\dim\M=n<\infty$, with parametrisation $\theta:\M\ra\Theta\subseteq\RR^n$, and constraints given by 
\[
	F(q)=\lambda_1\left(\int\dd x\int\dd\theta q(x|\theta)-1\right)+\lambda_2\left(\int\dd\theta q(x|\theta)-\delta(x-b)\right),
\]
where $\lambda_1$ and $\lambda_2$ are Lagrange multipliers.

\section{Towards new approach}  

In order to bypass the above problems, we need to take more careful look at the foundations of bayesian approach. Note that the Ramsey--de~Finetti type \cite{Ramsey:1931,deFinetti:1931,deFinetti:1937} and Cox's type \cite{Cox:1946,Cox:1961} derivations of Bayes' theorem (or, equivalently, of the algebraic rules of `probability calculus') \textit{assume} that the conditional probabilities $p(A|I)$ are to be used in order to draw inferences on the base of premises (evidence) $I$. Hence, the function $p(A|I)$, or any other function used to derive it, already \textit{assumes} that \textit{some} rule of probability updating (inductive inference) has to be used, because only this assumption allows to speak of some elements of the algebra as `evidence', or to speak of conditional probabilities as `inferences'. Any particular algebraic rules of transformation of conditional probabilities arise only as a result of \textit{additional} assumptions, which might not be relevant for general purposes and require anyway some additional justification. This observation allows us to consider spaces $\M(\mho)\subseteq L_1(\mho)^+$ of unconditioned positive finite integrals (\textit{information models}) as a candidate for kinematic component of inductive inference theory, and to consider \textit{some} principle $\PPP:\M(\mho)\ra\M(\mho)$ of updating of integrals as a candidate for its dynamic component (\textit{information dynamics}). As opposed to approaches aimed at identification of algebraic and lattice theoretic relations underlying evaluations $p(A|I)$ and their transformations \cite{Knuth:2009,Knuth:Skilling:2010}, we do not require that inferences should be conditioned exclusively on the elements of the underlying algebra, and we allow infinite-dimensional algebras and information models in foundations.

The choice of any particular form of principle of inductive inference (information dynamics) is a delicate issue, because (for any particular form of information kinematics) it determines the range and form of allowed inferences. According to the arguments of \cite{Hume:1748,Chwistek:1948}, there can be given no deductive logical premises for the claim that some inductive inference rule is absolute (universal), while inductive justification of induction is impossible due to circularity. However, if the chosen principle reduces in particular cases to a wide class of practically convenient and in some sense optimal techniques, then it can be considered as \textit{appealing}. Moreover, if it allows characterisation by axioms possessing unambiguous interpretation, then it can be considered as \textit{appealing} too. 

Because Bayes' rule is a special case of constrained maximisation of relative entropy \cite{Caticha:Giffin:2006}, we consider the latter as a candidate for a general principle $\PPP$ of quantitative inductive inference. For an evidence that it is appealing from the practical point of view, let us note that: (i) the conditional expectations are characterised as maximisers of the expectations of Bregman relative entropy \cite{BGW:2005}; (ii) the maximum likelihood methods are just special case of application of Bayes' rule \cite{Jaynes:2003}; (iii) inference techniques based on Fisher--Rao information metric amount to using the second order Taylor expansion of Csisz\'{a}r--Morimoto relative entropies \cite{Kullback:1959,Eguchi:1983}; (iv) many standard frequentist techniques of statistical inference can be reexpressed in terms of relative entropy, see e.g. \cite{Kullback:1959,Cressie:Read:1984,Zhu:Rohwer:1998,Jaynes:2003,Eguchi:Copas:2006}. Regarding axiomatisation, Shore and Johnson \cite{Shore:Johnson:1980,Shore:Johnson:1981}, Paris and Vencovsk\'{a} \cite{Paris:Vencovska:1989,Paris:Vencovska:1990}, and Csisz\'{a}r \cite{Csiszar:1991} have provided characterisations of the principle of maximisation of constrained Kullback--Leibler relative entropy as a unique probability updating rule that satisfies some set of conditions (see also \cite{Skilling:1988,Uffink:1995,Uffink:1996,Caticha:2004,Caticha:2007}). More generally, the existence and uniqueness of constrained maximisation of any  Bregman relative entropy is equivalent with the projection along $\nabla$-geodesic onto the non-empty set that is convex and closed in $\nabla^\nsdual$-geodesics, where the torsion-free flat affine connections $(\nabla,\nabla^\nsdual)$ arise from the third order Taylor expansion of this relative entropy \cite{Amari:1985,Amari:Nagaoka:1993}. If any of conditions of this type are accepted (what forms a particular \textit{decision}), then the resulting updating rule is unique, what makes this rule axiomatically appealing. However, like in the case of derivation of Bayes' rule from Cox's type  or the Ramsey--de Finetti type  procedure, one might deny some of the premises of these derivations (such as normalisation), and \textit{decide} to accept some other set of premises, leading to some other inductive inference rule. There also remains a problem how to relate the mathematical constraints of relative entropic updating with the operational descriptions specifying particular experiment. 

The choice of a particular form of information dynamics is thus relativised to a particular set of decisions, which are \textit{in principle} arbitrary. The same applies also to the choice of particular form of information kinematics (which includes model construction and model selection). But this arbitrariness is not necessarily unconstrained. According to the `subjective' bayesian interpretation, it is constrained by consistency of decisions of a single person (individual). Thus, each person can in principle choose arbitrary method of kinematic model construction and arbitrary method of inductive inference, but he or she is required to maintain personal consistency of these choices in subsequent inferences. In our opinion, the \textit{necessary} requirement for \textit{scientific} inference (as opposed to \textit{personal} inference) is to make these decisions consistent relatively to a particular community of users/agents. In other words, the decisions underlying construction and use of a particular form of information kinematics and information dynamics should be intersubjectively accepted and applied by all members of the given community. This way, within the range of intersubjective validity of these decisions, the notion of information model and its dynamics cannot be considered as `subjective'. 

This asks for the \textit{sufficient} conditions that define the \textit{scientific} inference. The crucial observation comes from Fleck \cite{Fleck:1935} (see also \cite{Spengler:1918} and \cite{Kuhn:1962}), who showed that the `scientific facts' and `experimental data' are always specified within the frames of some decisions (expressed in terms of particular assumptions and settings, including specification of the allowed response scales of measured outcomes, allowed experimental preparations and treatments, etc.), which are necessary to obtain intersubjective consistency with the preconceived notion of an `experiment of a given type'. These decisions \textit{define} the range of allowed variability of `facts' and `data', but are not determined by the theoretical model under consideration. Yet, \textit{both} theoretical (inferential, informational) and experimental (operational, knowledge-carrying) layers of scientific inquiry, as well as their mutual relationships, are determined by some intersubjectively shared \textit{thought style} (which is an entity from meta-theoretic level). Hence, `scientific facts' or `experimental data' are relative to some particular intersubjectively shared decisions on construction and use of experimental setups. Everything that is individually (\textit{personally}) experienced in a particular experimental situation, but does not fit into the frames rendered by the above decisions, is not considered as a valid `experimental data' (`scientific facts') for an experiment \textit{of a given type}. However, it \textit{does not} mean that the `experimental data' for a particular instance of experiment of a given type is completely \textit{determined} by these frames. Taking under consideration all above restrictions and relativisations, there remains unexpectedness of a particular outcome that appears \textit{within given frames} (of particular configuration and particular range of allowed outcomes). The aim of theoretical inquiry is to provide inductive inferences about these unexpected outcomes, dependently on the particular constraints that are taken into consideration as evidence. 

We postulate that for any (stage of historical development of a) given thought style, the scientific character of inquiry requires to separate  the theoretical abstract language used to intersubjectively define and communicate theories from the operational language used to intersubjectively define and communicate corresponding experiments. The specification of operational (`experimental') language for intersubjectively valid construction and use of experimental setups allows to define the range of validity of particular methods of inductive inference. More precisely, this allows to construct the \textit{relational} justification of given rules of construction of information kinematics and information dynamics with respect to a particular description of experiments of a given type, provided in terms of a given operational language. In particular, this applies to justification of the choice of $E$ and $F$ in the constrained maximisation of relative entropy. They can be determined by the assumptions expressed in operational language \textit{and} by the chosen relationship (correspondence) between operational and theoretical descriptions. The latter interprets terms of operational language in theoretical terms, and provides an operational (experimental) contextualisation of a given information (inductive inference) theory. 

This means that the theoretical model is verified only \textit{with respect to} certain context of intersubjectively shared decisions which construct the `experiment of a given type'. More precisely, the `experimental verification' of a predictive theory means just an \textit{intersubjective reproducibility (consistency)} of relationship between predictions (inferences) over a particular information model and results of use of an experimental setup of a given type, \textit{with respect to} a given correspondence between the kinematics of this model and construction of experimental setup of a given type, as well as between the constraints of inductive inference and the particular use of this experimental setup. To simplify these relativisations, we can restrict the notion of \textit{intersubjective validity}, as refering to fixed operational definitions \textit{and} fixed operational contextualisation. This leads to \textit{intersubjective interpretation}: the meaning of knowledge used to define particular theoretical model and its dynamics is provided by operational criteria that are sufficient and necessary in order to \textit{intersubjectively reproduce} an `experiment of a given type' that is considered to correspond to this theoretical model (which means that the inferences drawn from this model are interpreted as most plausible outcomes of corresponding experiment). 

This interpretation does not define the absolute (passive, static) meaning of the notion of `knowledge'. It defines only the relational (active, dynamic) meaning of this notion, as a particular relationship between kinematics-and-dynamics of theoretical model and construction-and-use of experimental setup. It differs from Bridgman's operationalism \cite{Bridgman:1927} and various versions of logical positivism, because it postulates neither the \textit{reduction} of theoretical layer to experimental layer, nor the absolute universality of their relationship. It differs from conventionalism of Duhem \cite{Duhem:1894,Duhem:1906}, Poincar\'{e} \cite{Poincare:1902,Poincare:1908}, and late Jeffreys \cite{Jeffreys:1955} by an additional requirement of intersubjective consistency between experimental setup construction and theoretical model construction. This way it is capable of providing a solution to the problem that neither `subjective' bayesianism nor `objective' bayesianism can justify the particular use of `experimental data' as evidence in inductive inference procedures. Personal betting behaviour and ontological postulates play no role in it. The intersubjective consistency (validity) of a particular relationship between theoretical model construction and operational construction of experimental setup, as well as a particular relationship between theoretical dynamics and operational use of experimental setup, is relative only to some community of users/agents which agree upon them. The constraint of intersubjective agreement is of meta-theoretical character and cannot be described in terms of inductive inference theory. Beyond any given community, the particular rules of models construction and inductive inference, as well as their relationship with the particular experimental setups and their use, are \textit{irrelevant} (arbitrary, personalistic, `subjective'), but within this range they are \textit{indispensable} (necessary, scientific, `objective'). This resolves the ``subjective vs objective'' bayesian debate, as well as it dissolves the bayesian version of the reference class problem \cite{Hajek:2007}, by removing them \textit{beyond} inductive inference theory. The criteria for intersubjective consistency (coherence) of experimental verifiability establish the direct link of inductive inference with experimental data, which is independent of any ontological assumptions.

\section{The principles of new foundations} 

On the level of \textit{mathematical framework}, we propose to unify kinematic (probabilistic, evaluational) and dynamic (statistic-inferential, relational) components, taking the best insights from analysis of the Borel--Kolmogorov and the Bayes--Laplace approaches. Thus, we replace measure spaces $(\X,\mho(\X),\mu)$ by camDcb-algebras $\mho$, and we consider integrals instead of measures. The failure of frequentism allows us to introduce statistical inference and lack of normalisation directly into foundations. We define: \textit{information kinematics} as given by \textit{information models} $\M(\mho)\subseteq L_1(\mho)^+$ and their \textit{information geometry} (specified by deviations, riemannian metrics, affine connections, etc., see e.g. \cite{Amari:Nagaoka:1993}); \textit{information dynamics} as given by minimisation of information deviation $D:\M(\mho)\times\M(\mho)\ra[0,\infty]$ such that $D(\omega,\phi)=0\iff\omega=\phi$, weighted by relative prior measure $E:\M(\mho)\times\M(\mho)\ra[0,\infty]$, and subjected to constraints $F:\M(\mho)\ra]-\infty,\infty]$,
\begin{equation}
	\PPP_{D,E,F}:\M(\mho)\ni\omega\mapsto\arginf_{\phi\in\M(\mho)}\left\{\int_{\varphi\in\M(\mho)} E(\varphi,\omega)D(\varphi,\phi)+F(\phi)\right\}\in\M(\mho).
\label{ent.upd}
\end{equation}

On the level of \textit{information semantics}, the underlying algebra $\mho$ represents an \textit{abstract qualitative language} subjected to quantitative evaluation, the space $\M(\mho)$ of finite integrals and its geometry represents \textit{quantified knowledge}, while the entropic updating \eqref{ent.upd} represents \textit{quantitative inductive inference procedure} determined by the triple $(D,E,F)$. The functions $E$ and $F$ specify the \textit{evidence}, while the resulting projection $\PPP_{D,E,F}$ is an \textit{inductive inference}: the most plausible state of knowledge subjected to given evidence. By \textit{well-posed inference problem} we understand such triples $(D,E,F)$ for which the solution $\PPP_{D,E,F}$ exists and is unique. The lack of existence indicates \textit{overdetermination} of an inference problem, while the lack of uniqueness indicates its \textit{undetermination}. For $E(\varphi,\omega)=\dd\varphi\delta(\varphi-\omega)$ the principle \eqref{ent.upd} says: given initial information state $\omega$, choose such information state that is most close to $\omega$ in terms of distance defined by $D$, under constraints defined by $F$. Other relative prior measures $E$ on $\M(\mho)$ allow more general selection of information state, which takes under consideration the relative distances to several different information states associated  with the initial state. 

To determine the particular operational and conceptual meaning attributed to the terms `knowledge' and `change of knowledge', the \textit{information semantics} requires an additional \textit{epistemic interpretation}. It should determine the choice of a particular information kinematics ($\M(\mho)$ and its information geometry) and a particular information dynamics ($D$, $E$, $F$) when applied to some particular experiments. According to \textit{intersubjective epistemic interpretation}: (i) the construction of information kinematics should correspond to the intersubjective description of construction of experimental setup of a given type (provided in some operational terms); (ii) the construction of information dynamics should correspond to the intersubjective description of the use of a experimental setup of a given type. But what do `intersubjective description' and `correspondence' mean? 

While it is impossible to formalise the intersubjective \textit{decisions} within the framework of inductive \textit{inference} theory, one can formalise the operational criteria that are necessary and sufficient in order to verify whether individual instances of experimentation (e.g., some mutually related preparations, actions and observations) can be considered as an \textit{intersubjectively valid} instance of an `experiment of a given type'. These criteria \textit{define} an intersubjective notion of an `experiment of a given type', which is an abstraction, because one cannot step twice into the \textit{same} stream. The intersubjective reproducibility (validation) of experiments can be undestood as a preservation of their abstract structure in all individual (personal) instances of experimentation. We will encode this structure together with the corresponding structure-preservation criteria in terms of category theory. By an `experiment of a given type' we understand an `experimental setup of a given type' \textit{and} a `given type of its use'. The former is defined by providing some particular categories: (i) $\Subj$ of abstract qualities (`preparations', `subjects', `things', `properties', `questions', `experimental units') subjected to consideration in experiment; (ii) $\Config$ of active configurations of experimental inputs (`treatments', `interventions', `acts'); (iii) $\Scale$ of registration scales of experimental outcomes (`outputs', `results'). The latter is defined by providing an association between inputs and outcomes conditioned on the choice of qualities. The particular association of outcome, input and quality in course of an individual (personal) instance of an experiment will be called `fact' or `experimental data'. Following \cite{Hinkelmann:Kempthorne:2005} we define an \textit{observation} as a map from preparation (quality) to outcome (registration value), and define an \textit{experimental design} as a map from preparation to treatment (act, configuration). This allows us to define \textit{fact} as a map from experimental design to an outcome. We also require that facts can be decomposed into `active' parts provided by treatments, and `passive' parts provided by observations. Following \cite{McCullagh:2002} we assume that $\Subj$, $\Config$ and $\Scale$ are equipped with forgetful functors $U_\cdot:\cdot\ra\Set$ (allowing to deal with their objects as with sets), and we define: a \textit{category of observations} as a product category $\Scale\times\Subj^{\mathrm{op}}$; a category of \textit{experimental designs} as a comma category $U_\Subj\downarrow U_\Config$; a \textit{category of facts} as 
\[
	\Scale\times(U_\Subj\downarrow U_\Config)^{\mathrm{op}}.
\]
The decomposion of facts into experimental designs, treatments and observations is provided, respectively, by means of the projection functors 
\begin{align*}
\mathbb{P}_{ed}&:=(\cdot)^{\mathrm{op}}\circ\Pi_{(U_\Subj\downarrow U_\Config)^{\mathrm{op}}},\\
\mathbb{P}_{t}&:=\Pi_{\Config},\\
\mathbb{P}_{o}&:=\id_\Scale\times\Pi_{\Subj^{\mathrm{op}}}.
\end{align*}

Let theoretical layer be equipped not only with models $\M(\mho)$ and triples $(D,E,F)$, but also with spaces $H$ of hypotheses. Define an \textit{inductive inference theory} as a method that assigns to each model $\M(\mho)$ a $(D,E,F)$-dependent map into some $H$. If this method can be expressed as a covariant functor from the category of models $\M(\mho)$ equipped with information geometric structures and geometry-preserving morphisms as arrows to the category of spaces of hypotheses, then by \textit{correspondence} we undestand the functorial association of the respective categories:
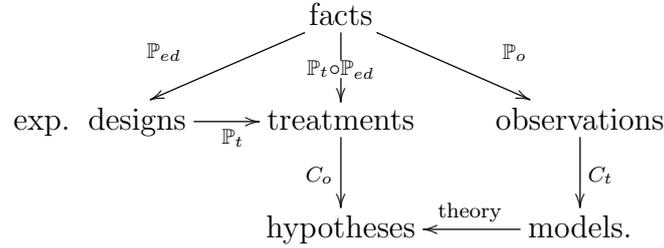
\[
       \xymatrix{
                &
                \mathrm{facts}
                \ar[ld]_{\mathbb{P}_{ed}\;\;\;\;\;\;\;\;\;\;}
                \ar[rd]^{\;\;\;\;\;\;\;\;\;\;\mathbb{P}_{o}}
                \ar[d]|{\mathbb{P}_t\circ\mathbb{P}_{ed}}
                &
                \\
                \mathrm{exp.\;\;designs}
                \ar[r]_{\;\;\;\;\mathbb{P}_{t}}
                &
                \mathrm{treatments}
                \ar[d]_{C_o}
                &
                \mathrm{observations}
                \ar[d]^{C_t}
                \\
                &
                \mathrm{hypotheses}
                &
                \mathrm{models}.
                \ar[l]_{\;\;\;\mathrm{theory}}
        }
\]
This scheme is inspired and motivated by \cite{Hinkelmann:Kempthorne:2005} and \cite{McCullagh:2002}. The functors $C_t$ and $C_o$ are  correspondence rules that interpret active and passive components of  experimental layer in terms of, respectively, active and passive components of theoretical layer. The theory will be called \textit{experimentally verifiable} with respect to the given intersubjective context of facts and correspondence rules iff the above diagram commutes. This replaces the `absolute' evaluation of the truth (verifiability, falsifiability) of a theory by the \textit{mutual coherence} of experimental and theoretical layers of intersubjective scientific inquiry in face of the given correspondence rules. In analogy to \cite{McCullagh:2002} one can consider the algebras $\mho$ as objects of $\Subj$. In such case each $\mho$ can be understood as an abstract entity that is used for intersubjective communication of experimental preparations. The preparations represented by $\mho$ are subjected to observations valued into the objects of $\Scale$. The interpretation $C_o$ maps these valuations into the corresponding information models $\M(\mho)$ and their geometry. The models $\M(\mho)$ consists of valuations (integrals) into $\RR^+$, and are understood as the carriers of quantitative intersubjective knowledge describing the experimental setup of a given type. The choice of evidence $E$ and $F$ provides the description of active configurational settings (treatments, inputs). The choice of $D$ and $\PPP_{D,E,F}$ determines a range of allowed inferences. The experimental verifiability condition determines the relationship of these inferences with experimental facts.

\nocite{Kostecki:2010:AIP}

{\small
\bibliographystyle{../rpkbib}
\bibliography{../rpkrefs}
}

\end{document}